\documentclass[preprint,review,12pt]{elsarticle}
\usepackage{epsfig} % for postscript graphics files
\usepackage{amsmath} % assumes amsmath package installed
\usepackage{amssymb}  % assumes amsmath package installed
\usepackage{ulem}
\usepackage{url}
\usepackage{amsmath,times,amssymb}
\usepackage{psfrag,graphicx}
\usepackage{subfigure}

\newtheorem{alg}{Algorithm}

\makeatletter
\renewcommand*\env@matrix[1][*\c@MaxMatrixCols c]{%
  \hskip -\arraycolsep
  \let\@ifnextchar\new@ifnextchar
  \array{#1}}
\makeatother

\allowdisplaybreaks

\label{newcom}
\newcommand{\xb}{\mathbf{x}}

\newcommand{\Hb}{\mathbf{H}}
\newcommand{\Pb}{\mathbf{P}}
\newcommand{\pb}{\mathbf{p}}

\newcommand{\rr}{\mathbb{R}}

\newcommand{\Nr}{\mathcal{N}}

\newcommand{\grad}{\nabla}

\begin{document}

\begin{frontmatter}

\title{\LARGE \bf
A distributed accelerated gradient algorithm for distributed model predictive control of a hydro power valley
}

\author{{\bf{Minh Dang Doan}}$^{*}$ {\bf{Pontus Giselsson}}$^{**}$ {\bf{Tam\'{a}s
  Keviczky}}$^{***}$\\
{\bf{Bart De Schutter}}$^{***}$ {\bf{Anders Rantzer}}$^{**}$}

\address{$^{*}$\phantom{a}{\it{Department of Electrical, Electronic, and Telecommunications Engineering, Cantho University of Technology, Vietnam}}\\
{\it{(email: minhdang@doan.vn)}}\\
$^{**}$\phantom{a}{\it{Department of Automatic Control, Lund University, Sweden}}\\
{\it{(email: \{pontusg,rantzer\}@control.lth.se})}\\
$^{***}$\phantom{a}{\it{Delft Center for Systems and Control, Delft University of Technology, The Netherlands}}\\
{\it{(email: \{t.keviczky,b.deschutter\}@tudelft.nl})}\\
}

%\maketitle
%\thispagestyle{empty}
%\pagestyle{empty}

%%%%%%%%%%%%%%%%%%%%%%%%%%%%%%%%%%%%%%%%%%%%%%%%%%%%%%%%%%%%%%%%%%%%%%%%%%%%%%%%
\begin{abstract}

A distributed model predictive control (DMPC) approach based on distributed optimization is applied to the power reference tracking problem of a hydro power valley (HPV) system. The applied optimization algorithm is based on accelerated gradient methods and achieves a convergence rate of $O\left(\frac{1}{k^2}\right)$, where $k$ is the iteration number. Major challenges in the control of the HPV include a nonlinear and large-scale model, nonsmoothness in the power-production functions, and a globally coupled cost function that prevents distributed schemes to be applied directly. We propose a linearization and approximation approach that accommodates the proposed the DMPC framework and provides very similar performance compared to a centralized solution in simulations. The provided numerical studies also suggest that for the sparsely interconnected system at hand, the distributed algorithm we propose is faster than a centralized state-of-the-art solver such as CPLEX.

\end{abstract}


\end{frontmatter}
%%%%%%%%%%%%%%%%%%%%%%%%%%%%%%%%%%%%%%%%%%%%%%%%%%%%%%%%%%%%%%%%%%%%%%%%%%%%%%%%
\section{Introduction}\label{sec_intro}

Hydro power plants generate electricity from potential energy and
kinetic energy of natural water, and often a number of power plants
are placed along a long river or a water body system to generate the
power at different stages. Currently, hydro power is one of the most
important means of renewable power generation in the world \citep{WB:2011_hydro}. In order to meet the world's electricity demand, hydro power
production should continue to  
grow due to the increasing cost of fossil fuels. However, hydro electricity, 
like any renewable energy, depends on the availability of a primary
resource, in this case: water. %Most natural locations where
%power-generating infrastructure can be built economically, have
%already been utilized \citep{PEW:11hydro}. 
The expected trend for
future use of hydro power is to build small-scale plants that can
generate electricity for a single community. Thus, an increasingly
important objective of hydro power plants is to manage the available
water resources efficiently, while following an optimal production
profile with respect to changes in the electricity market, to maximize
the long-term benefit of the plant. 
This water resource management must be compatible with ship navigation 
and irrigation, and it must respect environmental and safety
constraints on levels and flow rates in the lakes and the rivers. By significantly increasing the power efficiency of hydro power valley (HPV) systems, real-time control of water flows becomes an important ingredient in achieving this objective.

An HPV may contain several rivers and lakes, spanning a wide
geographical area and exhibiting complex dynamics. In order to tackle
the plant-wide control of such a complex system, an HPV is often
treated as a large-scale system consisting of interacting subsystems.
Large-scale system control has been an active research area that has
resulted in a variety of control techniques, which can be classified
in three main categories: decentralized control, distributed
control, and centralized control. The application of these approaches can be found
in a rich literature on control of water canals for irrigation and
hydro systems \citep{MarWey:05ARC, LitFro:09hdro}. We are interested
in applying model predictive control (MPC), a control method that has
been successfully used in industry \citep{QinBad:03CEP}, thanks to its
capability of handling hard constraints and the simple way of
incorporating an economical objective by means of an optimization
problem. For the control problem of open water systems, centralized
MPC has been studied in numerical examples using nonlinear MPC approaches in
combination with model smoothing and/or model reduction techniques
\citep{IgrLem:09LNCIS, NedSch:11NSC}, and in real implementations with
linear MPC of low-dimensional systems \citep{Overloop_thesis,
  OveCle:10JIDE}. However, centralized MPC has a drawback when
controlling large-scale systems due to limitations in communications
and the computational burden. These issues fostered the studies of
decentralized MPC and distributed MPC for large-scale water systems.
Early decentralized MPC methods for irrigation canals used the
decomposition-coordination approach to obtain decentralized versions
of LQ control \citep{FawGeo:98SMC}. Several decentralized MPC
simulations applied to irrigation canals and rivers were presented in
\citep{Georges:94CCA, Sawadogo:1998, GomRod:02AMM, SahMor:10}.
Distributed MPC approaches based on coordination and cooperation for
water delivery canals were presented in \cite{Georges:94CCA,
  NegOve:08_NHM, IgrCad:11MCCA, AnaJos:11NSC}. The typical control
objective in these studies is to regulate water levels and to deliver
the required amount of water to the right place at some time in the
future, i.e., the cost function does not have any special term except
the quadratic penalties on the states and the inputs. On the other
hand, in hydro power control, there are output penalty terms in the
cost function that represent the objective of manipulating power
production. Recent literature taking into account this cost function
includes centralized nonlinear MPC with a parallel version of the
multiple-shooting method for the optimal control problem using
continuous nonlinear dynamics \citep{SavRom:11JPC}, and a software
framework that formulates a discrete-time linear MPC controller with
the possibility to integrate a nonlinear prediction model and to use
commercial solvers to solve the optimization problem
\citep{Petrone:10MScThesis}. The hydro power control problem
considered in the current paper is similar to the setup in
\cite{SavRom:11JPC, Petrone:10MScThesis}. However, it distinguishes
itself by using a distributed control structure that aims to avoid
global communications and that divides the computational tasks into
local sub-tasks that are handled by subsystems, making the approach
more suitable for scaling up to even more complicated hydro power
plants.

The distributed MPC design approach proposed in this paper is enabled
by a distributed optimization algorithm that has
recently been developed by the authors in \cite{GisDoa:11Aut}. This
optimization algorithm is designed for a class of strongly convex
problems with mixed 1-norm and 2-norm terms in the cost function,
which perfectly suits the power reference tracking objective in the
HPV control benchmark. The underlying optimization algorithm in
\cite{GisDoa:11Aut}, although being implemented in a distributed way,
is proved to achieve the global optimum with an $O(\frac{1}{k^2})$
convergence rate, where $k$ is the iteration 
number. This is a significant
improvement compared to the distributed MPC methods presented in 
\cite{DoaKev:11JPC,DoaKev:09,giselssonDMPC,Negenborn07}, which achieve an 
$O(\frac{1}{k})$ convergence rate. There are three main challenges in
applying distributed MPC using the 
algorithm from \cite{GisDoa:11Aut} to the HPV benchmark problem. The
first one is that the nonlinear continuous-time model yields a
relatively large linear model after spatial and temporal
discretizations. We present a decentralized model order reduction
method that significantly reduces the model complexity while
maintaining prominent dynamics. The second challenge is that the
power production 
functions are nonsmooth, which prevents gradient-based methods to be
applied directly.
A method to overcome this difficulty and to enable optimal control using the
algorithm from \cite{GisDoa:11Aut} is also 
presented. The third challenge is that the whole
system should follow a centralized 
power reference which, if the algorithm from \cite{GisDoa:11Aut} is
applied directly, requires centralized communication. We propose a
dynamic power division approach that allows to track this centralized
power 
reference with only distributed communications. By means of numerical
examples, we will demonstrate the fast convergence 
property of the distributed algorithm which, when implemented on a
single core, can outperform a
state-of-the-art centralized solver (CPLEX) when solving the same
optimization problem.

%- Outline
The remaining parts of the paper are organized as follows. In
Section~\ref{sec_problem}, we describe the HPV system
and the power reference tracking problem that were formulated in the
HPV benchmark problem \citep{SavDie:11_hpv}.
Section~\ref{sec_dist_opt} provides a summary of the distributed
optimization framework that the authors have developed in
\cite{GisDoa:11Aut}. In Section~\ref{sec_control}, we present our
approach for modeling and model reduction of the HPV system, followed
by a reformulation of the MPC optimization problem, and developing a
distributed estimator so that the closed loop distributed MPC scheme
can be implemented using neighbor-to-neighbor communications only. The
simulation results are presented in Section~\ref{sec_simulations},
which also features a comparison with centralized MPC and
decentralized MPC. Through the various aspects of the comparison 
including performance, computational efficiency, and communication
requirements, the advantages of the distributed MPC algorithm will be
highlighted. Section~\ref{sec_conclusions} concludes the paper and
outlines future work.

\section{Problem description}\label{sec_problem}
In this section, we provide a summary of the hydro power valley
benchmark \citep{SavDie:11_hpv} and we present the linearized model
that serves as the starting point of our controller design.

\subsection{Hydro power valley system}
\begin{figure}[h!]
 \begin{center}
  \includegraphics[width=0.8\columnwidth]{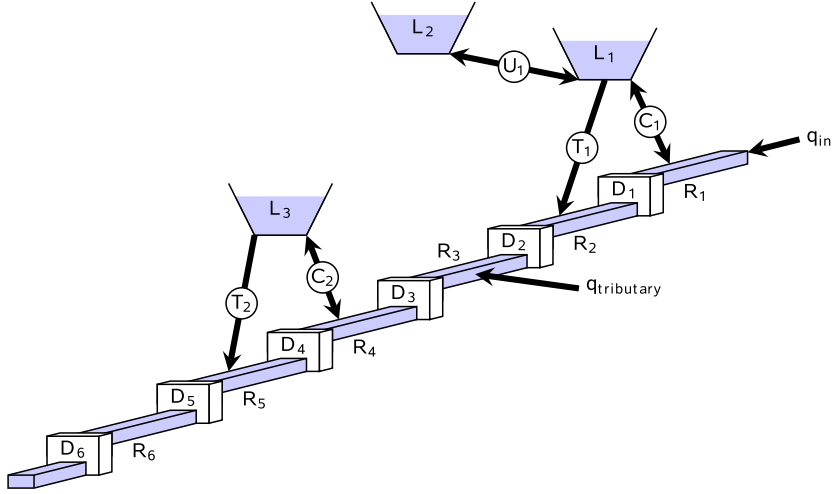}
  \caption{Overview of the HD-MPC hydro power valley system \citep{SavDie:11_hpv}}
  \label{fig_hpv}
 \end{center}
\end{figure}

We consider a hydro power plant composed of several interconnected
subsystems, as illustrated in Figure~\ref{fig_hpv}. The plant can be
divided into 8 subsystems, of which subsystem $S_1$ is composed of the
lakes $L_1, L_2$, the duct $U_1$ connecting them, and the ducts $C_1,
T_1$ that connect $L_1$ with the reaches\footnote{A reach is a river
  segment between two dams.} $R_1$, $R_2$, respectively. Subsystem
$S_2$ is composed of the lake $L_3$ and the ducts $C_2, T_2$ that
connect $L_3$ to the reaches $R_4, R_5$, respectively. There are 6
other subsystems, each of which consists of a reach and a dam at the
end of the reach. These six reaches $R_1, R_2, R_3, R_4, R_5$, and
$R_6$ are connected in series, separated by the dams $D_1, D_2, D_3,
D_4$, and $D_5$. The large lake that follows the dam $D_6$ is assumed
to have a fixed water level, which will absorb all the discharge. The
outside water flows enter the system at the upstream end of reach
$R_1$ and at the middle of reach $R_3$.

There are structures placed in the ducts and at the dams to control
the flows. These are the turbines placed in the ducts $T_1, T_2$ and
at each dam for power production. In the ducts $C_1, C_2$ there are
composite structures that can either function as pumps (for
transporting water to the lakes) or as turbines (when water is drained
from the lakes). 

The whole system has 10 manipulated variables, which are composed of
six dam flows ($q_{D1}$, $q_{D2}$, $q_{D3}$, $q_{D4}$, $q_{D5}$,
$q_{D6}$), two turbine flows ($q_{T1}$, $q_{T2}$) and two pump/turbine
flows ($q_{C1}$, $q_{C2}$). Further, the system has 9 measured
variables, the water levels in the three lakes ($h_{L1}$, $h_{L2}$,
$h_{L3}$) and the water levels at the end of each reach ($h_{R1}$,
$h_{R2}$, $h_{R3}$, $h_{R4}$, $h_{R5}$, $h_{R6}$).

\subsection{Power reference tracking problem}

One of the control problems specified in \cite{SavDie:11_hpv} is the
power reference tracking problem. We introduce state variables $x$,
which consist of water levels in the lakes and reaches and water
flows within the reaches, and control variables $q$, which are the
manipulated water flows.
The problem is to track a power production profile, $p^{\mathrm{ref}}(t)$, on a daily
basis using the following cost function:
\begin{align}\label{eq_tracking_cost}
J \triangleq &\int_{0}^{T} \gamma \left| p^{\mathrm{ref}}(t) -
  \sum_{i=1}^8 p_i(x(t),q(t))\right|dt \nonumber\\
&+ \sum_{i=1}^8 \int_{0}^{T} (x_i(t)-x_i^{\mathrm{ss}})^T Q_{i}
(x_i(t)-x_i^{\mathrm{ss}})dt \nonumber\\
&+ \sum_{i=1}^8 \int_{0}^{T} (q_i(t)-q_i^{\mathrm{ss}})^T R_i
(q_i(t)-q_i^{\mathrm{ss}})dt
\end{align}
subject to the nonlinear dynamics and linear constraints on
outputs and inputs as specified in
\cite{SavDie:11_hpv}. The
weights $Q_i, R_i, i=1,\dots,8$, $\gamma$, and the testing period $T$
are parameters of the benchmark.

The quadratic term in the cost function represents the penalties on the state deviation from the
steady state $x^{\mathrm{ss}}$ and the energy used for manipulating
the inputs away from the steady state flows $q^{\mathrm{ss}}$. The 1-norm term represents the power reference tracking mismatch, in which the function $p^{\rm{ref}}$ is the power reference and the
function $p_i$ represents the locally produced/consumed power by a subsystem
$i\in\{1,\dots,8\}$. For $i=1,2$ the produced/consumed power is (cf.\
\cite{SavDie:11_hpv})
\begin{equation}\label{eq:prodPowIeq1to2}
p_i(x(t),q(t))=k_{C_i}(q_{C_{i}}(t))q_{C_{i}}(t)\Delta
x_{C_i}(t)+k_{T_i}q_{T_{i}}(t)\Delta x_{T_i}(t)
\end{equation}
where $q_{C_i}$ and $q_{T_i}$ are the flows through ducts $C_i$ and $T_i$,
$\Delta x_{C_i}$ and $\Delta x_{T_i}$ are the relative differences in water
levels before and after ducts $C_i$ and $T_i$ respectively, $k_{T_i}$ is the power coefficient of the turbine $T_i$, and 
\begin{align}\label{eq_k_C1C2}
k_{C_i}(q_{C_{i}}(t))=\left\{\begin{array}{ll}
k_{T_{C_i}},~ & q_{C_i}(t)\geq 0\\
k_{P_{C_i}},~ & q_{C_i}(t)<0\\
\end{array}\right.
\end{align}
is a discontinuous power coefficient that depends on whether the duct $C_i$ acts as a turbine ($q_{C_i}(t) \geq 0$) or as a pump ($q_{C_i}(t)<0$). For $i=3,\ldots,8$ we have
\begin{equation}\label{eq:prodPowIeq3to8}
p_i(x(t),q(t))=k_{D_{i-2}}q_{D_{i-2}}(t)\Delta x_{D_{i-2}}(t)
\end{equation}
which is the power produced by the turbine located at dam $D_{i-2}$. The
produced/consumed power functions given in \eqref{eq:prodPowIeq1to2} and \eqref{eq:prodPowIeq3to8} are nonlinear, and even nonsmooth for subsystems 1 and 2 due to the differences of $k_{T_{C_i}}$ and $k_{P_{C_i}}$ in \eqref{eq_k_C1C2}, thus complicating a direct application of a standard MPC scheme.

Still, the complexity of the system and control objective suggests an
optimization based control strategy, such as MPC.
Further, the distributed nature of the system makes it possible to consider
distributed MPC techniques. However, the stated optimization problem
\eqref{eq_tracking_cost} is a nonlinear continuous-time dynamic
optimization problem, which in general is very hard to solve. In the next sections we will discuss the modeling of the hydro power valley that leads to a linearized model.
%In the following sections we show how to approximate this optimization problem to fit the distributed MPC formulation presented in \cite{GisDoa:11Aut}. 

\subsection{Nonlinear hydro power valley model}\label{HPVnonlinModel}
The model of the reaches is based on the one-dimensional Saint Venant
partial differential equation, representing the mass and momentum
balance (see \cite{SavDie:11_hpv} for details): 
\begin{equation}\label{eq:sve}
\left\{
\begin{array}{l}
\dfrac{\partial q(t,z)}{\partial z} + \dfrac{\partial s(t,z)}{\partial t} = 0 \\[0.2cm]
\dfrac{1}{g}\dfrac{\partial}{\partial t}\left(\dfrac{q(t,z)}{s(t,z)}\right) +
\dfrac{1}{2g}\dfrac{\partial}{\partial z}\left(\dfrac{q^2(t,z)}{s^2(t,z)}\right) +
\dfrac{\partial h(t,z)}{\partial z} + I_\mathrm{f}(t,z) - I_0(z) = 0
\end{array}
\right.
\end{equation}
with $z$ the spatial variable, $t$ the time variable, $q$ the river
flow (or discharge), $s$ the cross-section surface of the river, $h$
the water level w.r.t. the river bed, $I_\mathrm{f}$ the friction
slope, $I_0(z)$ the river bed slope, and $g$ the gravitational
acceleration constant.

%Spatial Discretization 
The partial differential equation \eqref{eq:sve} is converted into a
system of ordinary differential equations by using spatial
discretization. To achieve this, each reach is divided into 20 cells,
yielding 20 additional states, which are the water levels at the
beginning of the cells. For details of the spatial discretization and the
equations for the resulting nonlinear dynamical system the reader is referred to
\cite[Section 2.1.1]{SavDie:11_hpv}. The resulting nonlinear dynamical
system has in total 249 states, 10 inputs, and 9 outputs.

\subsection{Model linearization and discretization}

As mentioned in Section~\ref{HPVnonlinModel} a set of nonlinear ordinary differential equations that describe
the hydro power valley dynamics is presented in \cite[Section
2.1.1]{SavDie:11_hpv}. A linear continuous-time model which is
linearized around the steady
state operating point $(x^{\mathrm{ss}}, q^{\mathrm{ss}})$ is also
provided in the HPV benchmark package \citep{SavDie:11_hpv}. 
Discretizing this model using zero-order-hold gives a
discrete-time linear system with 249 states and 10
inputs. The coupling of the subsystems is through the inputs only. This
implies that discretization using
zero-order-hold of the continuous-time system keeps the structure of the
original system description. Thus, the resulting discrete time system
has a block-diagonal dynamics matrix, a block-diagonal output matrix, and
a sparse input matrix, and each subsystem 
$i=1,\dots,8$ can be expressed in the following form:
\begin{align}\label{eq_lin_model_local}
x_i^{\mathrm{d}}(k+1)&= A_{ii} x_i^{\mathrm{d}}(k)+\sum_{j=1}^{8} B_{ij}q_j^{\mathrm{d}}(k) \\
y_i^{\mathrm{d}}(k) &= C_{i}x_i^{\mathrm{d}}(k) \nonumber
\end{align}
in which the variables $x^{\mathrm{d}}, q^{\mathrm{d}}$, and $y^{\mathrm{d}}$ stand for the deviation
from the steady-state values, and the subscripts $i,j$ stand for the
subsystem indices. As mentioned the subsystems are coupled through the
inputs only and at least for some $j\in\{1,\ldots,8\}$ we have
$B_{ij}=0$ for every $i=1,\ldots,8$.

The use of a discrete-time linearized model
enables controller design with some
specific approaches, which include our proposed distributed
optimization technique presented in \cite{GisDoa:11Aut}. Before
describing our main contributions, we now provide a summary of this
distributed optimization framework in the next section.

\section{Distributed optimization framework for MPC}\label{sec_dist_opt}

In this section, we describe the distributed optimization algorithm developed in
\cite{GisDoa:11Aut} which is based on an accelerated gradient method.
The first accelerated gradient method was developed in 
\cite{Nesterov1983} and further elaborated and extended in 
\cite{BecTab_FISTA:2009,Nesterov1988,Nes_smooth:2005,TohYun_acc:2010,Tseng_acc:2008}.
The main idea of the algorithm presented in \cite{GisDoa:11Aut} is to
exploit the problem structure of the
dual problem such that accelerated gradient computations can be distributed to
subsystems. Hence, the distributed algorithm effectively solves the centralized optimization problem. 
Dual decomposition has been used in the past to tackle the complexity of large-scale optimization problems arising in water supply networks \citep{carpentier+93}. In our work however, in addition to simplifying the local computations, we apply this decomposition philosophy in order to distribute the decision-making process.

The algorithm in \cite{GisDoa:11Aut} is developed to handle
optimization problems of the form
\begin{align}\label{eq:optProb}
\min_{\xb,\xb_{\mathrm{a}}} ~& \frac{1}{2}\xb^T\Hb\xb+g^T\xb+\gamma\|\xb_{\mathrm{a}}\|_1\\
\textrm{s.t.} ~& \mathbf{A}\xb=\mathbf{b} \nonumber\\
& \mathbf{C}\xb\leq \mathbf{d} \nonumber\\
&\xb_{\mathrm{a}}=\Pb\xb-\pb\nonumber
\end{align}
where $\xb\in\mathbb{R}^n$ and $\xb_{\mathrm{a}}\in\mathbb{R}^m$ are
vectors of decision variables,
and $\xb$ is
partitioned according to:
\begin{align}
 \xb = [\xb_1^T, \dots, \xb_M^T]^T,
\label{eq:xPart}
\end{align}
and $\xb_i\in\mathbb{R}^{n_i}$. 
Further, the matrix $\Hb\in\mathbb{R}^{n\times n}$ is positive
definite and block-diagonal, the matrices $\mathbf{A} \in \rr^{q\times
  n}$, $\mathbf{C} \in
\rr^{r\times n}$, and $\Pb \in \rr^{m\times n}$ have sparse
structures, and $g\in\mathbb{R}^n$, $\pb\in\mathbb{R}^m$,
$\mathbf{b}\in\mathbb{R}^q$, $\mathbf{d}\in\mathbb{R}^r$.
We introduce the partitions $g = [g_1^T,\ldots,g_M^T]^T$, $\pb =
[\pb_1^T,\ldots,\pb_M^T]^T$, $\mathbf{b} =
[\mathbf{b}_1^T,\ldots,\mathbf{b}_M^T]^T$, $\mathbf{d} =
[\mathbf{d}_1^T,\ldots,\mathbf{d}_M^T]^T$,
\begin{align*}
\Hb &= \begin{bmatrix}
\Hb_1 & &\\
& \ddots &\\
& & \Hb_M
\end{bmatrix},& \mathbf{A} &= \begin{bmatrix}
\mathbf{A}_{11} & \ldots & \mathbf{A}_{1M}\\
\vdots & \ddots & \vdots\\
\mathbf{A}_{M1} & \ldots & \mathbf{A}_{MM}
\end{bmatrix}\\
\mathbf{C} &= \begin{bmatrix}
\mathbf{C}_{11} & \ldots & \mathbf{C}_{1M}\\
\vdots & \ddots & \vdots\\
\mathbf{C}_{M1} & \ldots & \mathbf{C}_{MM}
\end{bmatrix},& \Pb &= \begin{bmatrix}
\Pb_{11} & \ldots & \Pb_{1M}\\
\vdots & \ddots & \vdots\\
\Pb_{M1} & \ldots & \Pb_{MM}
\end{bmatrix}
\end{align*}
where the partitions are introduced in accordance with
\eqref{eq:xPart} and $g_i\in\mathbb{R}^{n_i}$,
$\pb_i\in\mathbb{R}^{m_i}$, $\mathbf{b}_i\in\mathbb{R}^{q_i}$,
$\mathbf{d}_i\in\mathbb{R}^{r_i}$, $H_i\in\mathbb{R}^{n_i\times n_i}$,
$\mathbf{A}_{ij}\in\mathbb{R}^{q_i\times n_j}$,
$\mathbf{C}_{ij}\in\mathbb{R}^{r_i\times n_j}$ and
$\Pb_{ij}\in\mathbb{R}^{m_i\times n_j}$.
The assumption on sparsity of $\mathbf{A}$,
$\mathbf{C}$ and $\Pb$ is that $\mathbf{A}_{ij}=0$,
$\mathbf{C}_{ij}=0$, and $\Pb_{ij}=0$ for some $i, j$ and we assume that
the constraint matrices are built such that $\mathbf{A}_{ii}\neq 0$,
$\mathbf{C}_{ii}\neq 0$, and $\Pb_{ii}\neq 0$ for all $i=1,\ldots,M$.
Based on the coupling, we define for each subsystem a neighborhood
set, denoted by $\Nr_i$, as follows: 
\begin{align}\label{eq_Ni_def}
  \Nr_i = \big\{j \in \{1, \dots, M\} |~ \mathbf{A}_{ij} \neq 0 {\hbox{ or }}
 \mathbf{A}_{ji}\neq 0 {\hbox{ or }} \mathbf{C}_{ij} \neq 0 {\hbox{ or }}
 &\mathbf{C}_{ji}\neq 0 {\hbox{ or }} \\
\nonumber & \Pb_{ij} \neq 0 {\hbox{ or }}
 \Pb_{ji}\neq 0 \big\}.
\end{align}
Note that there are two type of equality constriants in \eqref{eq:optProb}, the first one involves only $\xb$ and the matrix $\mathbf{A}$ has a sparsity pattern, i.e., there is no global coupling introduced in that equality constraint; the last one involves both $\xb$ and $\xb_{\mathrm{a}}$, moreover introduces a global coupling due to the fact that $\xb_{\mathrm{a}}$ is penalized in the 1-norm term of the cost function, thus it is not straightforward to deal with this constraint as we could treat the first constraint. Throughout the paper, the dual variables corresponding to these constraints are treated differently, and a distributed approximation of the 1-norm term is introduced to treat the second type of equality constraint.

We introduce dual variables $\lambda\in\mathbb{R}^q,
\mu\in\mathbb{R}^r, \nu\in\mathbb{R}^m$ for the equality
constraints, inequality constraints, and equality constraints
originating from the 1-norm cost in \eqref{eq:optProb} respectively.
We also introduce the dual variable partitions $\lambda =
[\lambda_1^T,\ldots,\lambda_M^T]^T$, $\mu =
[\mu_1^T,\ldots,\mu_M^T]^T$, and $\nu = [\nu_1^T,\ldots,\nu_M^T]^T$
where $\lambda_i\in\mathbb{R}^{q_i}$, $\mu_i\in\mathbb{R}^{r_i}$, and
$\nu_i\in\mathbb{R}^{m_i}$. Based on \cite{GisDoa:11Aut}, the dual problem of
\eqref{eq:optProb} can be cast as the minimization of the negative dual function
\begin{align}
\label{eq:convDualFcn}f(\lambda,\mu,\nu) 
=\frac{1}{2}(\mathbf{A}^T\lambda+\mathbf{C}^T\mu+P^T\nu)^T\Hb^{-1}(\mathbf{A}^T\lambda+\mathbf{C}^T\mu&+\Pb^T\nu)+\\
\nonumber &+\mathbf{b}^T\lambda+\mathbf{d}^T\mu+\pb^T\nu
\end{align}
and the dual variables are constrained to satisfy
\begin{align}
\lambda&\in\mathbb{R}^q, & \mu&\in\mathbb{R}_{+}^{r},&
\nu&\in[-\gamma,\gamma]^{m}
\end{align}
% \begin{equation}\label{eq:setZ}
% z \in Z = \mathbb{R}^{q}\times \mathbb{R}_{+}^{r} \times [-\gamma,\gamma]^{m}
% \end{equation}
where $\mathbb{R}_{+}$ denotes the non-negative real orthant. 
The negative dual function \eqref{eq:convDualFcn} has a Lipschitz continuous
gradient with constant (cf. \cite{GisDoa:11Aut})
\begin{align}\label{eq_Lipschitz_const}
 L=\|[\mathbf{A}^T~\mathbf{C}^T~P^T]^T\Hb^{-1}[\mathbf{A}^T~\mathbf{C}^T~\Pb^T]\|_2
\end{align}
and can hence be minimized using accelerated gradient methods. The
distributed accelerated gradient method as presented in
\cite{GisDoa:11Aut} is summarized below in a slightly different form that is adapted to our HPV application problem at hand.

\begin{alg}\label{alg:accGrad} \textbf{Distributed accelerated
    gradient algorithm}
\hrule
\vspace{3mm}
\noindent Initialize $\lambda^0=\lambda^{-1}$,
$\mu^0=\mu^{-1}$, $\nu^0=\nu^{-1}$ and $\xb^{-1}$ with the last values
from the previous sampling step. For the first sampling step, these
variables are initialized by zeros.\\
In every node, $i$, the following computations are performed:\\
{\bf{For}} $k=0,1,2,\dots$
\begin{enumerate}
 \item Compute
\begin{align*}
 \xb_i^{k} &=
 -\Hb_i^{-1}\bigg(\sum_{j\in\mathcal{N}_i}\left(\mathbf{A}_{ji}^T\lambda_j +\mathbf{C}_{ji}^T\mu_j+\Pb_{ji}^T\nu_j\right)\bigg)\\
\bar{\xb}_i^{k}&=\xb_i^k+\frac{k-1}{k+2}(\xb_i^k-\xb_i^{k-1})
\end{align*}
 \item Send $\bar{\xb}_i^{k}$ to each $j \in \mathcal{N}_i$, receive
   $\bar{\xb}_j^{k}$ from each $j \in \Nr_i$\\
 \item Compute
\begin{align*}
\lambda_i^{k+1} &= \lambda_i^k+\frac{k-1}{k+2}(\lambda_i^k-\lambda_i^{k-1})+\frac{1}{L}\bigg(\sum_{j\in\mathcal{N}_i}\mathbf{A}_{ij}\mathbf{\bar{x}}_j-\mathbf{b}_i\bigg)\\
\mu_i^{k+1} &=
\max\bigg\{0,\mu_i^k+\frac{k-1}{k+2}(\mu_i^k-\mu_i^{k-1})+\frac{1}{L}\bigg(\sum_{j\in\mathcal{N}_i}\mathbf{C}_{ij}\mathbf{\bar{x}}_j-\mathbf{d}_i\bigg)\bigg\} \\
\nu_i^{k+1} &= \min\bigg\{\gamma,\max\bigg[-\gamma,
\nu_i^k+\frac{k-1}{k+2}(\nu_i^k-\nu_i^{k-1})+\\
&\qquad\qquad\qquad\qquad\qquad\qquad\qquad\qquad+
\frac{1}{L}\bigg(\sum_{j\in\mathcal{N}_i} \Pb_{ij}\mathbf{\bar{x}}_j-\pb_i\bigg)\bigg]\bigg\} 
\end{align*}
 \item Send $\lambda_{i}^{k+1}$, $\mu_{i}^{k+1}$,
   $\nu_{i}^{k+1}$ to each $j \in \Nr_i$, 
   receive $\lambda_{j}^{k+1}$, $\mu_{j}^{k+1}$,
   $\nu_{j}^{k+1}$ from each $j\in
   \mathcal{N}_i$.
\end{enumerate}
\end{alg}
\hrule
\bigskip

The Lipschitz constant $L$ of $\grad f$ is used
in the algorithm. For MPC purposes we only need to compute $L$
once in a centralized way and use it through all MPC problem
instances.

Besides the suitability for distributed implementation, another merit
of Algorithm~\ref{alg:accGrad} is its fast convergence rate. The main
convergence results of Algorithm~\ref{alg:accGrad} are given in
\cite{GisDoa:11Aut}, stating that both the dual function value and the primal
variables converge towards their respective optima with the rate of
$O\left(\frac{1}{k^2}\right)$ where $k$ is the iteration index. This
convergence rate is much better than the convergence rate of classical
gradient-based optimization algorithms, which is
$O\left(\frac{1}{k}\right)$.

\section{Control of HPV using distributed MPC}\label{sec_control}

We have so far described the linear discrete-time model of the HPV
in Section~\ref{sec_problem} and the fast distributed optimization
method, Algorithm~\ref{alg:accGrad}, that serves as a basis for
designing a distributed model
predictive controller to be applied to the HPV. However, there are
three major challenges for this application. First, the linear
discrete-time model cannot be directly used in an MPC context due to
the existence of a number of 
unobservable and uncontrollable modes. These
unobservable/uncontrollable modes are a result of the spatial discretization in each reach which creates states that cannot be observed/controlled separately. In addition, the linear discrete-time model has a large
number of states, causing a large computational burden. Second, the power functions associated 
with the ducts $C_1$ and $C_2$ are nonsmooth (cf.
\eqref{eq:prodPowIeq1to2} and \eqref{eq_k_C1C2}). The nonsmoothness is
caused by the fact that the flow through $C_1$ and $C_2$ is bidirectional and 
the powers consumed/produced do not have equivalent coefficients. The
third major challenge is the global coupling in the cost function due
to the fact that we have to track a central power reference function that specifies the desired sum of locally generated power outputs. This
global coupling prevents a distributed implementation of
Algorithm~\ref{alg:accGrad} since the sparsity in the constraints is
lost. These issues are addressed in the following sections.

\subsection{Modification of the linear model}\label{sec_modeling}
In this section we show how to create a model of the HPV that is suitable for the
DMPC framework presented in \cite{GisDoa:11Aut}. First we present a
model reduction technique that keeps the system structure, then the
nonsmooth power function is treated.

\subsubsection{Decentralized model order reduction}

The block-diagonal structure of discrete-time dynamical system
\eqref{eq_lin_model_local} makes it possible to perform model reduction on
each subsystem individually. Several model reduction methods have been proposed for interconnected systems \citep{VanDoo:2007chap, SanMur:2009OCAM}. In this work, we use a straightforward balanced
truncation method \citep{GugAnt:04IJC,Moore:81TAC} to reduce the order of each local model
\eqref{eq_lin_model_local}.

Let us introduce $B_{i} = [B_{i1} \dots
B_{i8}]$ and $q^{\mathrm{d}} = [(q_1^{\mathrm{d}})^T \dots (q_8^{\mathrm{d}})^T]^T$ to get the following
discrete-time linear model of each subsystem:
\begin{align}\label{eq_subsystem_model}
x_i^{\mathrm{d}}(k+1)&= A_{ii}x_i^{\mathrm{d}}(k)+B_{i}q^{\mathrm{d}}(k) \\
y_i^{\mathrm{d}}(k) &= C_{i}x_i^{\mathrm{d}}(k).\nonumber
\end{align}

Applying the balanced truncation technique yields transformation
matrices denoted by $T_i^{\mathrm{r}}$ and
$T_i^{\mathrm{r},\mathrm{inv}}$ for each subsystem, where
$T_i^{\mathrm{r}} T_i^{\mathrm{r},\mathrm{inv}} = I$. By denoting the
new state variables, $x_i^{\mathrm{r}}=T_i^rx_i^{\mathrm{d}}$, and the
control variable $q^{\mathrm{r}}=q^{\mathrm{d}}$, we represent the reduced order
model as:
\begin{align}
x_i^{\mathrm{r}}(k+1)&= A_{ii}^{\mathrm{r}}x_i^{\mathrm{r}}(k)+B_{i}^{\mathrm{r}}q^{\mathrm{r}}(k) \label{eq_subsystem_model_reduced}\\
y_i^{\mathrm{r}}(k) &=
C_{i}^{\mathrm{r}}x_i^{\mathrm{r}}(k)\label{eq_subsystem_model_reduced_output}
\end{align}
where $A_{ii}^{\mathrm{r}} = T_i^{\mathrm{r}}A_{ii}T_i^{\mathrm{r,inv}}$, $B_{i}^{\mathrm{r}} = T_i^{\mathrm{r}} B_i$
and $C_{i}^{\mathrm{r}} = C_iT_i^{\mathrm{r,inv}}$.
It should be noted that the block-sparsity structure of $B_i^r$ is the same
as in the non-reduced input matrix $B_i$, since the model reduction is
performed for each local model separately. Moreover, all the modes of the reduced model are both observable and controllable.

The model reduction gives a 32-state reduced model that approximately
represents the dynamics of the full linear model with 249 states.

\subsubsection{Treatment of nonlinear and nonsmooth power function}\label{sec_virtual_flows}

One of the difficulties in applying a linear MPC approach to the hydro
power valley is the nonsmoothness of the power function associated
with the ducts $C_1$ and $C_2$, which is included in
the expression for power generation \eqref{eq:prodPowIeq1to2} in
subsystem 1 and subsystem 2, respectively.
In order to handle this nonsmoothness, we use a double-flow technique, which means
introducing two nonnegative positive variables to express the flow in
$C_i, i=1, 2$  at a sampling step $k$:
\begin{itemize}
\item $q_{C_{i\mathrm{P}}}(k)$: virtual flow such that $C_i$ functions as a pump
\item $q_{C_{i\mathrm{T}}}(k)$: virtual flow such that $C_i$ functions as a turbine
\end{itemize}
The introduction of virtual flows requires the input-matrices,
$B_i^{\mathrm{r}}$, to be augmented with two extra columns identical to the ones multiplying
$q_{C_i}, i=1,2$ with the opposite sign to capture that pump action is
also introduced with a positive flow. The resulting reduced order model
has 12 inputs instead of the original 10.
Using the introduced flows $q_{C_{i\mathrm{P}}}$ and
$q_{C_{i\mathrm{T}}}$, the power function \eqref{eq:prodPowIeq1to2} for subsystems 1 and
2 can be rewritten as
\begin{align}\label{eq_power_subsys_1_2}
p_i(x(k),q(k))=\left(k_{T_{C_i}}q_{C_{i\mathrm{T}}}(k)-k_{P_{C_i}}q_{C_{i\mathrm{P}}}(k)\right)\Delta
x_{C_i}(k)+k_{T_i}q_{T_{i}}(k)\Delta
x_{T_i}(k)
\end{align}
with the additional constraints that $q_{C_{i\mathrm{T}}}(k)\geq 0, q_{C_{i\mathrm{P}}}(k)\geq 0$ and
$q_{C_{i\mathrm{T}}}(k)q_{C_{i\mathrm{P}}}(k)=0$. The last constraint expresses the fact that water
flows in only one direction at a time, i.e., that either the pump or the
turbine is active. The resulting nonlinear expression \eqref{eq_power_subsys_1_2} can in
turn be linearized around the steady-state solution
$(x^{\mathrm{ss}},q^{\mathrm{ss}})$. Since $q_{C_i}^{\mathrm{ss}}=0$
for $i=1,2$ we get the following linear local power
production/consumption approximation for subsystems $i=1,2$:
\begin{multline*}
\hat{p}_i(x(k),q(k))=\Delta
x_{C_i}^{\mathrm{ss}}\left[k_{T_{C_i}}\;-k_{P_{C_i}}\right]
\begin{bmatrix}
q_{C_{i\mathrm{T}}}(k)\\
q_{C_{i\mathrm{P}}}(k)
\end{bmatrix}+\\+
k_{T_i}q_{T_{i}}^{\mathrm{ss}}\left(\Delta x_{T_i}(k)-\Delta x_{T_i}^{\mathrm{ss}}\right)+
k_{T_i}\Delta x_{T_i}^{\mathrm{ss}}\left(q_{T_{i}}(k)-q_{T_{i}}^{\mathrm{ss}}\right)+\\+k_{T_i}q_{T_i}^{\mathrm{ss}}\Delta x_{T_i}^{\mathrm{ss}}
\end{multline*}
This reformulation results in a linear expression with a nonlinear
constraint at each time step $k$, $q_{C_{i\mathrm{T}}}(k)q_{C_{i\mathrm{P}}}(k)=0$, that
approximates the original nonsmooth nonlinear power
production/consumption expression~\eqref{eq:prodPowIeq1to2}. We show
our approach to handle the nonlinear constraint in
Section~\ref{sec_optimization_formulation}.

For subsystems $i=3,\ldots,8$ we have smooth
power production expressions \eqref{eq:prodPowIeq3to8} that can be directly
linearized without introducing virtual flows:
\begin{multline*}
\hat{p}_i(x(k),q(k))=k_{D_i}q_{D_i}^{\mathrm{ss}}\Delta x_{D_i}^{\mathrm{ss}}+
k_{D_i}q_{D_{i}}^{\mathrm{ss}}\left(\Delta x_{D_i}(k)-\Delta
  x_{D_i}^{\mathrm{ss}}\right)+\\
+k_{D_i}\Delta x_{D_i}^{\mathrm{ss}}\left(q_{D_{i}}(k)-q_{D_{i}}^{\mathrm{ss}}\right)
\end{multline*}

\subsection{HPV optimization problem formulation}\label{sec_optimization_formulation}

In this section we formulate an optimization problem
of the form \eqref{eq:optProb} that can be used for power reference
tracking in the HPV benchmark using MPC. We have obtained a linear
discrete-time dynamical system
\eqref{eq_subsystem_model_reduced}-\eqref{eq_subsystem_model_reduced_output}
for the HPV with state variables
$x^{\mathrm{r}}$ and control variables $q^{\mathrm{r}}$. The
constraints are upper and lower bounds on the outputs and inputs and their values
can be found in \cite{SavDie:11_hpv}. Using the transformations matrices $T_i^{\mathrm{r}}$ and $T_i^{\mathrm{r},\mathrm{inv}}$, these constraints can readily be recast as linear constraints for the reduced order
model variables $x^{\mathrm{r}}, q^{\mathrm{r}}$.
The 
power reference problem formulation \eqref{eq_tracking_cost} specifies
a quadratic cost on states and control variables and a 1-norm
penalty on deviations
from the provided power reference, $p^{\mathrm{ref}}$. For
control horizon, $N$, this optimization problem can be
written as
\begin{align}\label{eq:HPVopt}
\min_{\mathbf{x},\mathbf{x}_{\mathrm{a}}} ~& \sum_{t=0}^{N-1}\left\{\sum_{i=1}^8\left[ x_i^{\mathrm{r}}(k)^T Q_i x_i^{\mathrm{r}}(k)+q_i^{\mathrm{r}}(k)^T R_i q_i^{\mathrm{r}}(k)\right]+\gamma\|x_{\mathrm{a}}(k)\|_1\right\}\\
\textrm{s.t.} ~& \begin{tabular}[t]{lll}
\eqref{eq_subsystem_model_reduced}, \eqref{eq_subsystem_model_reduced_output} & $k=0,\ldots,N-1$ & $i=1,\ldots,8$\\
$C^{\mathrm{r}}_i x_i^{\mathrm{r}}(k)\in \mathcal{Y}_i$ & $k=0,\ldots,N-1$ & $i=1,\ldots,8$\\
$q_i(k)\in \mathcal{Q}_i$ & $k=0,\ldots,N-1$ & $i=1,\ldots,8$\\
\multicolumn{2}{l}{$x_{\mathrm{a}}(k)=p^{\mathrm{ref}}(k)-\sum_{i=1}^8
  \hat{p}_i(x^{\mathrm{r}}(k),q^{\mathrm{r}}(k))$} & $k=0,\ldots,N-1$\\
$q_{C_{i\mathrm{T}}}(k)q_{C_{i\mathrm{P}}}(k)=0$ & $k=0,\ldots,N-1$ & $i=1,\ldots,2$
\end{tabular}\nonumber
\end{align}
where $\mathcal{Y}_i$ and $\mathcal{Q}_i$ are sets representing the
local output and input constraints, the additional variable $x_{\mathrm{a}}$ captures the power reference tracking mismatch, and the notation $\mathbf{x}$ represents the stack of variables $x_i^{\mathrm{r}}(k)$ and $q_i^{\mathrm{r}}(k)$ for all $i$ and $k$, while $\mathbf{x}_{\mathrm{a}}$ is the stacked variable of $x_{\mathrm{a}}(k)$ for all $k$. Note that we can write $\xb = [\xb_1^T, \dots, \xb_8^T]^T$ where each $\xb_i, i=1,\dots,8$ includes all the variables that belong to subsystem $i$.

\subsubsection{Power reference division}

Since the original cost function contains a non-separable 1-norm term,
the power reference constraints in the optimization problem
\eqref{eq:HPVopt} are coupled between all subsystems. This implies
that Algorithm~\ref{alg:accGrad} requires some 
global communication even though the only information needed to be
sent to the 
global coordinator is $\bar{p}_i(x^{\mathrm{r}}(k),q^{\mathrm{r}}(k))$
for $k=0,\ldots,N-1$ 
from each subsystem $i=1,\ldots,8$.

In order to obtain a suitable dual problem, we first need to reformulate the cost function in a separable form. For the sake of brevity, we focus on one sampling step and drop the time index $k$. Thus for now our simplified objective is to decompose the following problem:
\begin{align}\label{eq_cost_power}
 \min_{\{\xb_i\}_{i=1,\dots,8}} \quad \bigg|p^{\mathrm{ref}} - \sum_{i=1}^8 P_i \xb \bigg|
\end{align}
with $\xb = [\xb_1^T, \dots, \xb_8^T]^T$, and $P_i$ the matrix coefficient such that the power function produced or consumed by each subsystem 
$\hat{p}_i(x^{\mathrm{r}}(k),q^{\mathrm{r}}(k))$ is linearized as $P_i \xb(k)$.

In this section we present two different ways that avoid global
communication when solving this problem. In the first approach, we divide and distribute the global
power reference to the subsystems in a static fashion. In the second approach,
we show how the subsystems can trade local power references between
neighbors to achieve a satisfactory centralized reference tracking.

\paragraph{Static local power references}

The idea here is straightforward. We divide the global
power reference into local ones, i.e., $p^{\mathrm{ref}}$ is divided into
local parts $p_i^{\mathrm{ref}}$, $i=1,\ldots,8$. We have chosen to compute
$p_i^{\mathrm{ref}}$ such that it satisfies
\begin{align}
 \frac{p_i^\mathrm{ref}(k)}{\sum_{i=1}^8 p_i^\mathrm{ref}(k)}=\frac{p_i(x^{\mathrm{ss}},q^{\mathrm{ss}})}{\sum_{i=1}^8 p_i(x^{\mathrm{ss}},q^{\mathrm{ss}})}, \quad {\hbox{for }} k=0,\ldots,N-1
\end{align}
with $p_i(x^{\mathrm{ss}},q^{\mathrm{ss}})$ the power produced by subsystem $i$ in the steady-state condition.

This means that the fraction of the total power reference given to
subsystem $i$ is constant. The optimization problem is changed
accordingly, i.e., the following cost function can be used instead of \eqref{eq_cost_power}:
\begin{align}\label{eq_tracking_cost_group}
\min_{\{\xb_i\}_{i=1,\dots,8}} \quad \sum_{i=1}^{8} \bigg|p_i^{\mathrm{ref}} - P_i \xb \bigg|
\end{align}
with $\xb = [\xb_1^T, \dots, \xb_8^T]^T$. This allows for a distributed implementation since the matrix $P_i$ introduces only local couplings, i.e., subsystem $i$ needs only neighboring and local water levels and local water flows to
compute the corresponding power output. The disadvantage of the static power reference division is that the global
power reference tracking is not very accurate, as will be shown in the simulations section.

\paragraph{Dynamic local power references}

The static power division essentially means that each subsystem always
tracks a fraction of power reference that is equal to the proportion
it produces in the steady-state condition. When the total power
reference deviates significantly from the steady-state power, this
idea may not work well since the proportional change of the local
power reference can lead to sub-optimal performance. Inspired by an
idea in \cite{MadMar:11ACC}, we now introduce the dynamic power
division, in which the subsystems have more flexibility in choosing
the appropriate local power reference to be tracked. The main idea is
that each subsystem will exchange power references with
its direct neighbors.

Let us define for each pair $(i,j)$ with $j \in \Nr_i$ a node that is in charge of determining the power exchange variable between subsystems $i$ and $j$, denoted by $\delta_{ij}$ if node $i$ is in charge and by $\delta_{ji}$ if node $j$ is in charge \footnote{Note that here we discuss the power division for each sampling step, i.e., there are $\delta_{ij}(k)$ or $\delta_{ji}(k)$ with $k=0,\dots,N-1$.}. Then for each subsystem we form the set \footnote{A simple way is to let the subsystem with smaller index lead the exchange, i.e., $\Delta_i = \{j | j \in \Nr_i, j>i\}$.}:
\begin{align}
 \Delta_i = \{j~|~j \in \Nr_i, i \mathrm{~is~in~charge~of~} \delta_{ij}\}.
\end{align}
Now we replace \eqref{eq_cost_power} by the following cost function:
\begin{align}\label{eq_tracking_cost_exchange}
\min_{\{\xb_i,\mathbf{\delta}_i\}_{i=1,\dots,8}} \sum_{i=1}^{8} \bigg|p_i^{\mathrm{ref}} + \sum_{j \in \Delta_i} \delta_{ij} - \sum_{j \in \Nr_i \setminus \Delta_i} \delta_{ji} - P_i \xb \bigg|
\end{align}
with $\mathbf{\delta}_i$ the vector containing all $\delta_{ij}, j\in
\Delta_i$, and $p_i^{\mathrm{ref}}$ the nominal power reference for
subsystem $i$. In words, the local power reference for each subsystem
$i$ deviates from the nominal value by adding the exchange amounts of
the links that $i$ manages and subtracting the exchange amounts of the
links that affect $i$ but are decided upon by its neighbors. Note
that problem \eqref{eq_tracking_cost_exchange} has a sparse structure
that complies with the existing sparse structure of the HPV system,
i.e., this method does not expand the neighborhood set of each
subsystem.

The advantage of this dynamic power division is that it makes use of
the existing network topology to form a sparse cost function, and the
total power reference is preserved even if the local power references
can deviate from the nominal values, i.e., we always have:
\begin{align}
\sum_{i=1}^8 \bigg\{ p_i^{\mathrm{ref}} + \sum_{j \in \Delta_i} \delta_{ij} - \sum_{j \in \Nr_i \setminus \Delta_i} \delta_{ji} \bigg\} = p^{\mathrm{ref}}
\end{align}
\medskip

Now that we have a separable cost function by using either a static or
a dynamic power division technique, we can cast the approximate
optimization problem in the form \eqref{eq:optProb} that has a
separable dual problem, and apply Algorithm~\ref{alg:accGrad} at every
sampling step. However, due to the requirement of positive
definiteness of the quadratic term in the objective function, the
introduced power exchange variables $\delta_{ij}$ must be penalized
with a positive definite quadratic term. This implies that power
reference exchange has an associated cost.

\paragraph{Communication structures}

In the preceding sections we have presented three different ways to handle the
power reference term. The first is the one with centralized
power reference term which we hereby denote by GLOBAL--REF. The second is the
one with static local power references which we denote by LOC--REF--STAT.
The third is the dynamic local power reference which from here on
is denoted by
LOC--REF--DYN. In Table~\ref{tab_neighbors} we
provide an overview of the neighborhood sets $\mathcal{N}_i$ for the
different power reference tracking schemes.

\begin{table}[h!]
\centering
\caption{Neighborhoods of subsystems ($\mathcal{N}_i$)}
{\footnotesize{
\begin{tabular}{|c|c|c|c|}
\hline Subsystem & GLOBAL--REF &  LOC--REF--DYN & LOC--REF-STAT \\ 
\hline 1 & $\{1,\dots,8\}$ &  $\{1,3,4\}$ & $\{1,3,4\}$\\ 
%\hline 
2 & $\{1,\dots,8\}$ &  $\{2,6,7\}$ & $\{2,6,7\}$\\ 
%\hline 
3 & $\{1,\dots,8\}$ &  $\{3,1,4\}$ & $\{3,1,4\}$\\ 
%\hline 
4 & $\{1,\dots,8\}$ &  $\{4,1,3,5\}$ & $\{4,1,3,5\}$\\ 
%\hline 
5 & $\{1,\dots,8\}$ &  $\{5,4,6\}$ & $\{5,4,6\}$\\ 
%\hline 
6 & $\{1,\dots,8\}$ &  $\{6,2,7,5\}$ & $\{6,2,7,5\}$\\ 
%\hline 
7 & $\{1,\dots,8\}$ &  $\{7,2,6,8\}$ & $\{7,2,6,8\}$\\ 
%\hline 
8 & $\{1,\dots,8\}$ &  $\{8,7\}$ & $\{8,7\}$\\
\hline 
\end{tabular}}}
\label{tab_neighbors}
\end{table}
We can see that all subsystems have the same neighborhood sets for the dynamic local reference tracking and the static local
reference tracking.

\subsubsection{Relaxation of nonlinear constraint}\label{sec:NLconstr}

The second issue that hinders the optimization problem \eqref{eq:HPVopt}
from being solved using Algorithm~\ref{alg:accGrad} are the nonlinear
constraints $q_{C_{i\mathrm{T}}}(k)q_{C_{i\mathrm{P}}}(k)=0$ with $i=1, 2$. In this section we present a way
to relax these constraints.

Assuming in the cost function we have the penalty $R_{C_i} [q_{C_{i\mathrm{T}}} q_{C_{i\mathrm{P}}}]^T$ on the pump
and turbine action in ducts $C_i$, $i=1, 2$, with
\begin{equation}\label{eq:origCostNL}
R_{C_i} = \begin{bmatrix}
R_{C_{i\mathrm{T}}} & 0\\
0 & R_{C_{i\mathrm{P}}}
\end{bmatrix}.
\end{equation}
We also have the constraints that $q_{C_{i\mathrm{P}}}(k)\geq 0,
q_{C_{i\mathrm{T}}}(k)\geq
0$ and $q_{C_{i\mathrm{T}}}(k)q_{C_{i\mathrm{P}}}(k)=0$. We relax this by removing the
nonlinear constraint and adding a cross-penalty
$\alpha\sqrt{R_{C_{1\mathrm{P}}}R_{C_{1\mathrm{T}}}}$ for some $\alpha\in
(0,1)$ in the cost function, i.e., we set
\begin{equation}\label{eq:relaxedCostNL}
R_{C_i} = \begin{bmatrix}
R_{C_{i\mathrm{T}}} & \alpha\sqrt{R_{C_{i\mathrm{P}}}R_{C_{i\mathrm{T}}}}\\
\alpha\sqrt{R_{C_{i\mathrm{P}}}R_{C_{i\mathrm{T}}}} & R_{C_{i\mathrm{P}}}
\end{bmatrix}.
\end{equation}
This relaxation is implementable using the proposed algorithm since
the nonlinear constraint is removed and replaced by a cross-penalty.
The cross-penalty gives an additional
cost if both $q_{C_{i\mathrm{T}}}$ and $q_{C_{i\mathrm{P}}}$ are non-zero. The closer
$\alpha$ is to 1, the larger the penalty.
For $\alpha\geq 1$ it is easily verified that we lose strong
convexity on the quadratic cost function, i.e., $R_{C_i}$
loses positive definiteness and such choices for $\alpha$ are therefore
prohibited.

The relaxation is not equivalent to the original nonlinear constraint and thus
cannot guarantee that the nonlinear constraint is respected using
this relaxation. However, it
turns out that the optimal solution using
the cross-penalty in the cost \eqref{eq:relaxedCostNL} in most simulated cases
coincides with the optimal solution 
when the nonlinear constraint $q_{C_{i\mathrm{T}}}(k)q_{C_{i\mathrm{P}}}(k)=0$ and the
original diagonal cost \eqref{eq:origCostNL} are enforced. In some cases however, the optimal solution using the
relaxation does not respect the nonlinear constraint. To address
this, a two-phase optimization strategy is developed and presented next.

\subsubsection{Two-phase optimization}

We propose a two-phase optimization strategy as an ad-hoc branch and bound
optimization routine that uses two consecutive optimizations. In the
first optimization the relaxed optimization problem is solved. 
If the nonlinear constraints are
respected, i.e., we get a solution
that satisfies $q_{C_{i\mathrm{T}}}(k)q_{C_{i\mathrm{P}}}(k)=0$, the global optimal
solution for the non-relaxed problem is found.
If some of the nonlinear constraints do not hold, the optimization
routine is restarted with setting the smaller flow between $q_{C_{i\mathrm{T}}}(k)$ and
$q_{C_{i\mathrm{P}}}(k)$ to zero, for $i=1,2$, $k=0,\ldots,N-1$. The
resulting algorithm is summarized below.

\begin{alg}\label{BNB_alg} \textbf{Distributed branch and bound algorithm}
\hrule
%\vspace{3mm}
\begin{enumerate}
 \item Solve the relaxed problem using Algorithm~\ref{alg:accGrad}
 \item {\bf{If}} $q_{C_{i\mathrm{T}}}(k)q_{C_{i\mathrm{P}}}(k) \neq 0$, $\qquad\qquad i=1,2$, $t=0,\ldots,N-1$
   \begin{enumerate}
     \item[] {\bf{If}} $q_{C_{i\mathrm{T}}}(k) > q_{C_{i\mathrm{P}}}(k)$
       \begin{enumerate}
       \item[] Add constraint: $q_{C_{i\mathrm{P}}}(k)=0$
       \end{enumerate}
       \item[] {\bf{Else}}
         \begin{enumerate}
         \item[] Add constraint: $q_{C_{i\mathrm{T}}}(k)=0$
         \end{enumerate}
       \item[] {\bf{End}}
       \end{enumerate}
     \item[] {\bf{End}}
 \item Solve relaxed problem using Algorithm~\ref{alg:accGrad}
   with the additional flow constraints
\end{enumerate}
\end{alg}\hrule
\bigskip

This ad-hoc branch and bound technique does not theoretically guarantee that
the optimal flow directions are chosen. However, we can guarantee that
the nonlinear constraints are 
always satisfied. Further, for the distributed MPC
formulation we will see in the simulations section that the global
optimal solution for the non-relaxed problem is found at every time
step using this branch and bound algorithm.

\subsection{Distributed estimation}\label{sec:distrEst}

From Section \ref{sec_problem} we know that not all states can be measured, which
implies that an observer needs to be used to feed an initial condition
to the optimizer. The reduced-order linear model
\eqref{eq_subsystem_model_reduced}-\eqref{eq_subsystem_model_reduced_output}
has local
dynamics and outputs only, which implies that an observer can be
designed in decentralized fashion. We introduce the local estimate
$\hat{x}_i^{\mathrm{r}}$ and the local observer-gain $K_i$, and the
following local observer dynamics
\begin{equation*}
\hat{x}_i^{\mathrm{r}}(k+1) = A_{ii}^{\mathrm{r}}\hat{x}_i^{\mathrm{r}}(k)+B_i^{\mathrm{r}} q^{\mathrm{r}}(k)+K_i(y_i^{\mathrm{r}}(k)-C_i^{\mathrm{r}} \hat{x}_i^{\mathrm{r}}(k))
\end{equation*}
Because of the sparse structure of $B_i^{\mathrm{r}}$ this observer can be
implemented in a distributed fashion where only the inflows to subsystem $i$
need to be communicated. The estimation error $\widetilde{x}_i^{\mathrm{r}}=x_i^{\mathrm{r}}-\hat{x}_i^{\mathrm{r}}$ has
local error dynamics
\begin{equation*}
\widetilde{x}_i^{\mathrm{r}}(k+1) = (A_{ii}^{\mathrm{r}}-K_iC_i^{\mathrm{r}})\widetilde{x}_i^{\mathrm{r}}(k)
\end{equation*}
Thus, the observer can be designed in a decentralized fashion and be
implemented in a distributed fashion.

\section{Simulation results}\label{sec_simulations}

We perform distributed MPC simulations of the hydro power valley
using 3 different ways of handling the power reference: GLOBAL--REF,
LOC--REF--DYN, and 
LOC--REF--STAT, using the proposed Algorithm~\ref{BNB_alg}.
We also solve the problem \eqref{eq:HPVopt} using a state-of-the-art MIQP-solver, namely CPLEX. In CPLEX the nonlinear
constraints given in \eqref{eq:HPVopt} can be addressed by introducing binary variables. More specifically, for each duct $C_i, i=1, 2$, we define two virtual flows, $q_{C_{i\mathrm{P}}}$ and $q_{C_{i\mathrm{T}}}$, and require that both values are nonnegative. Each virtual flow has a maximum capacity, hence the constraints for these flows are:
\begin{align}\label{eq_const_virtual_flows}
 \left. \begin{array}{l}
0 \leq q_{C_{i\mathrm{P}}} \leq q_{C_{i\mathrm{P}}}^{\mathrm{max}} \\
0 \leq q_{C_{i\mathrm{T}}} \leq q_{C_{i\mathrm{T}}}^{\mathrm{max}}
\end{array}\right.
\end{align}
We introduce binary variables $b_i \in \{0,1\}$ and impose the following constraints:
\begin{align}\label{eq_const_bin_var}
 \left. \begin{array}{l}
q_{C_{i\mathrm{T}}} \leq q_{C_{i\mathrm{T}}}^{\mathrm{max}} b_i \\
q_{C_{i\mathrm{P}}} \leq q_{C_{i\mathrm{P}}}^{\mathrm{max}} (1-b_i)
\end{array}\right.
\end{align}

The constraints \eqref{eq_const_virtual_flows} and \eqref{eq_const_bin_var} ensure that either $q_{C_{i\mathrm{P}}}=0, q_{C_{i\mathrm{T}}} \geq 0$ (if $b_i=1$) or $q_{C_{i\mathrm{T}}}=0, q_{C_{i\mathrm{P}}} \geq 0$ (if $b_i=0$).

This formulation results in an MIQP for which there are efficient Branch-and-Bound algorithms
implemented in CPLEX. To make the 1-norm term in \eqref{eq:HPVopt} fit
the MIQP-formulation used in CPLEX we introduce auxiliary
variables $v$ and use the following equivalent reformulation
\begin{align*}
\min_x&\|Px-p\|_1 \qquad\Leftrightarrow &\min_{x,v}&~1^Tv\\
&&{\mathrm{s.t.}}& -v\leq Px-p\leq v
\end{align*}
We also compare the proposed distributed MPC method to a decentralized MPC approach in
which each subsystem solves its own local MPC problem without any
communication, in order to show the advantage of DMPC w.r.t.\ decentralized
MPC.

\subsection{Simulation details}

We use the original nonlinear continuous model presented in
\cite{SavDie:11_hpv} as simulation model. The ode-solver
\textit{ode15s} in MATLAB is used to perform the simulations. 
A MATLAB function that computes the derivatives needed by \textit{ode15s}
is provided in the benchmark package \citep{SavDie:11_hpv}.
The control system consists of the
distributed observer from Section~\ref{sec:distrEst} which feeds
Algorithm~\ref{BNB_alg}, with estimates of the current state.

Besides the mismatch
between the model used for control and the model used for simulation
we have also added 
bounded process noise to capture mismatch between the simulation
model and the real plant. The magnitude of the worst case process
noise was chosen to be $1\%$ of the steady-state level $x^{\mathrm{ss}}$.
We also use
bounded additive measurement noise where the measured water levels are
within $\pm 3$ cm from the actual water levels. 

We use a sampling time of 30 minutes in all simulations and the control horizon
is $N=10$, i.e., 5 hours. The simulations are
performed over a 24 hour period since the power reference trajectories
are periodic with this interval.

All simulations and optimizations were implemented on a PC running MATLAB on Linux
with an Intel(R) Core(TM) i7 CPU running at 3 GHz and with 4 GB RAM.

\subsection{Control performance comparison}

% \begin{figure}[htb]
% \centering
% \subfigure[Decentralized MPC]{
% %\includegraphics[width=0.45\columnwidth]{power_tracking_decMPC.eps}
% \includegraphics[width=0.45\columnwidth]{power_ref_decMPC_2.eps}
% \label{fig_tracking_decMPC}
% }
% \subfigure[DMPC and LOC--REF--STAT]{
% %\includegraphics[width=0.45\columnwidth]{power_tracking_DMPC_loc_stat_ref.eps}
% \includegraphics[width=0.45\columnwidth]{power_ref_loc_stat_2.eps}
% \label{fig_tracking_loc_stat}
% }
% \subfigure[DMPC and LOC--REF--DYN]{
% %\includegraphics[width=0.45\columnwidth]{power_tracking_DMPC_loc_dyn_ref.eps}
% \includegraphics[width=0.45\columnwidth]{power_ref_loc_dyn_2.eps}
% \label{fig_tracking_loc_dyn}
% }
% \subfigure[DMPC and GLOBAL--REF]{
% %\includegraphics[width=0.45\columnwidth]{power_tracking_DMPC_global_ref.eps}
% \includegraphics[width=0.45\columnwidth]{power_ref_central2.eps}
% \label{fig_tracking_global}
% }
% \label{fig_tracking}
% \caption{Comparison of power reference tracking performance using DMPC and decentralized MPC approaches.}
% \end{figure}

\begin{figure}[t]
\centering
\subfigure[Decentralized MPC]{
\psfrag{time}{{\footnotesize{$t$ [h]}}}
\psfrag{powerMWpow}{{\footnotesize{Power [MW]}}}
%\psfrag{aaaaaaaaaaaaaaaaaa}{{\scriptsize{Produced power}}}
%\psfrag{bbbbbbbbbbbbbbbbbb}{{\scriptsize{Reference power}}}
%\psfrag{cccccccccccccccccc}{{\scriptsize{Steady state power}}}
\includegraphics[width=5.3cm]{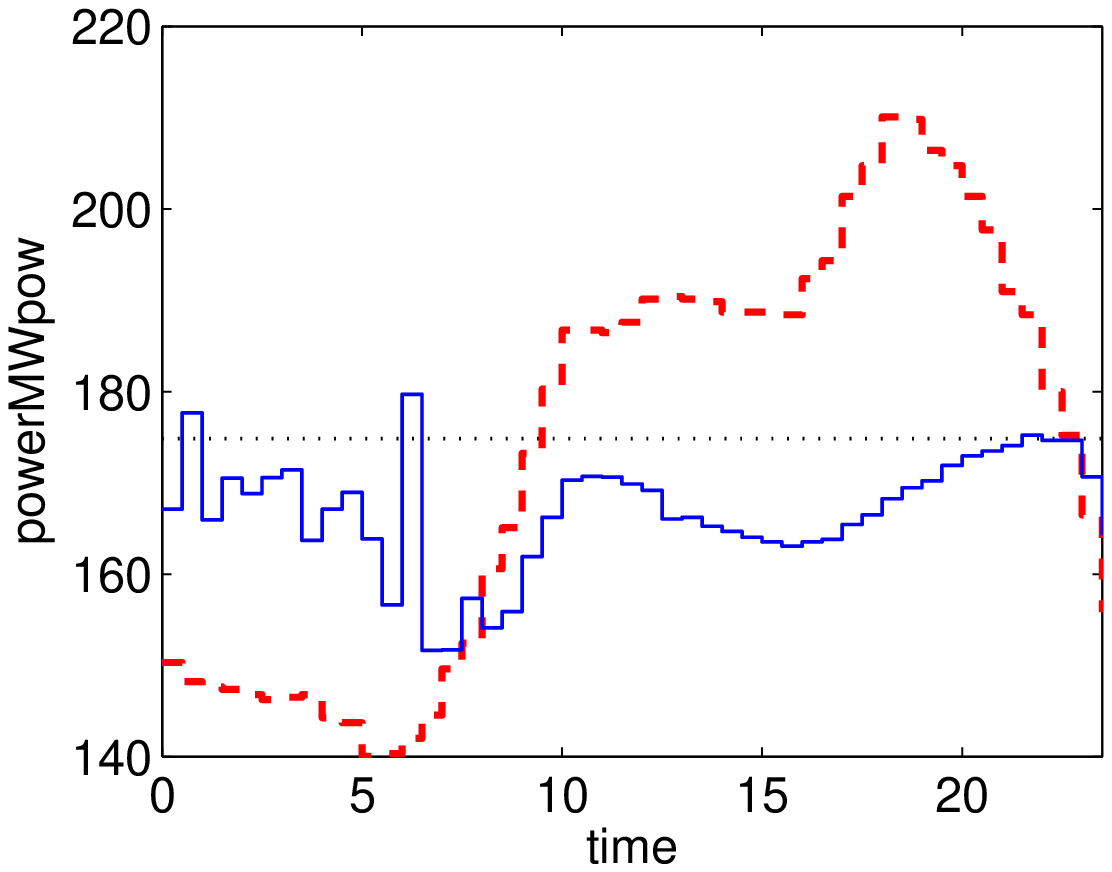}
\label{fig_tracking_decMPC}
}
\subfigure[DMPC and LOC--REF--STAT]{
\psfrag{time}{{\footnotesize{$t$ [h]}}}
\psfrag{powerMWpow}{{\footnotesize{Power [MW]}}}
%\psfrag{aaaaaaaaaaaaaaaaaa}{{\scriptsize{Produced power}}}
%\psfrag{bbbbbbbbbbbbbbbbbb}{{\scriptsize{Reference power}}}
%\psfrag{cccccccccccccccccc}{{\scriptsize{Steady state power}}}
\includegraphics[width=5.3cm]{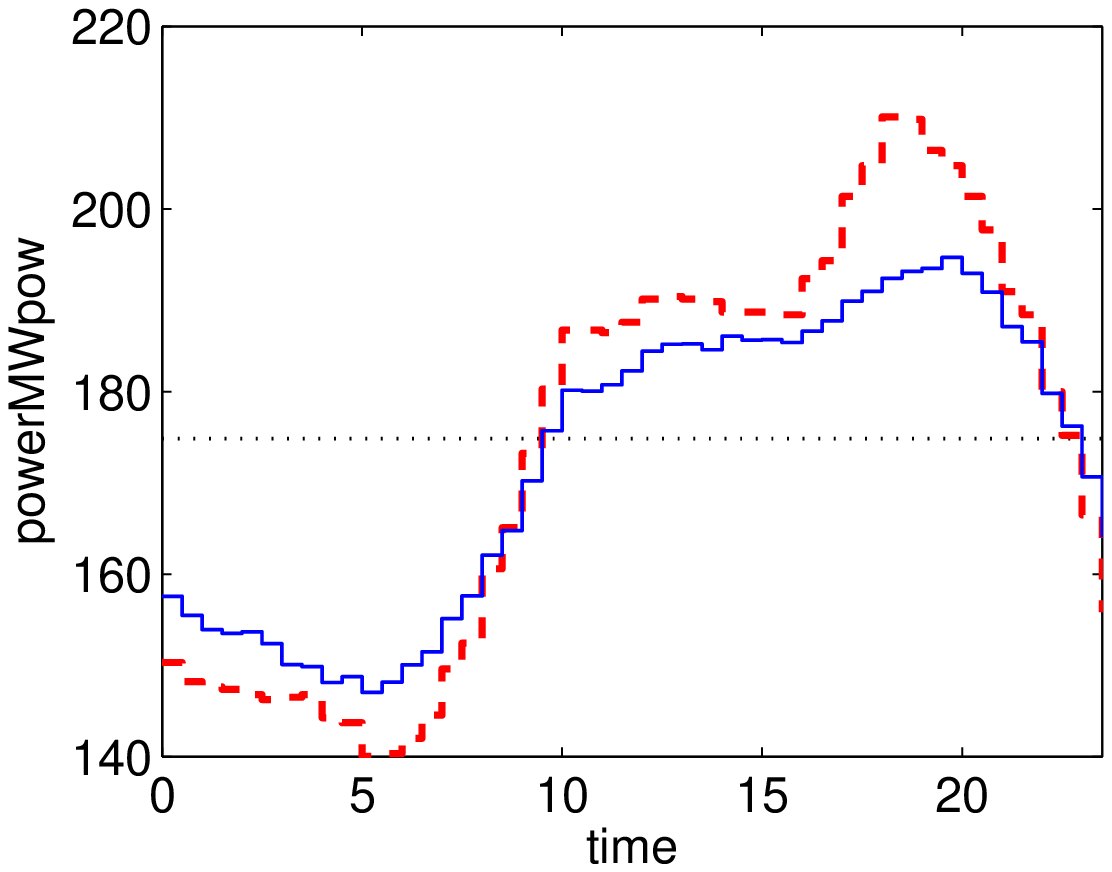}
\label{fig_tracking_loc_stat}
}
\subfigure[DMPC and LOC--REF--DYN]{
\psfrag{time}{{\footnotesize{$t$ [h]}}}
\psfrag{powerMWpow}{{\footnotesize{Power [MW]}}}
%\psfrag{aaaaaaaaaaaaaaaaaa}{{\scriptsize{Produced power}}}
%\psfrag{bbbbbbbbbbbbbbbbbb}{{\scriptsize{Reference power}}}
%\psfrag{cccccccccccccccccc}{{\scriptsize{Steady state power}}}
\includegraphics[width=5.3cm]{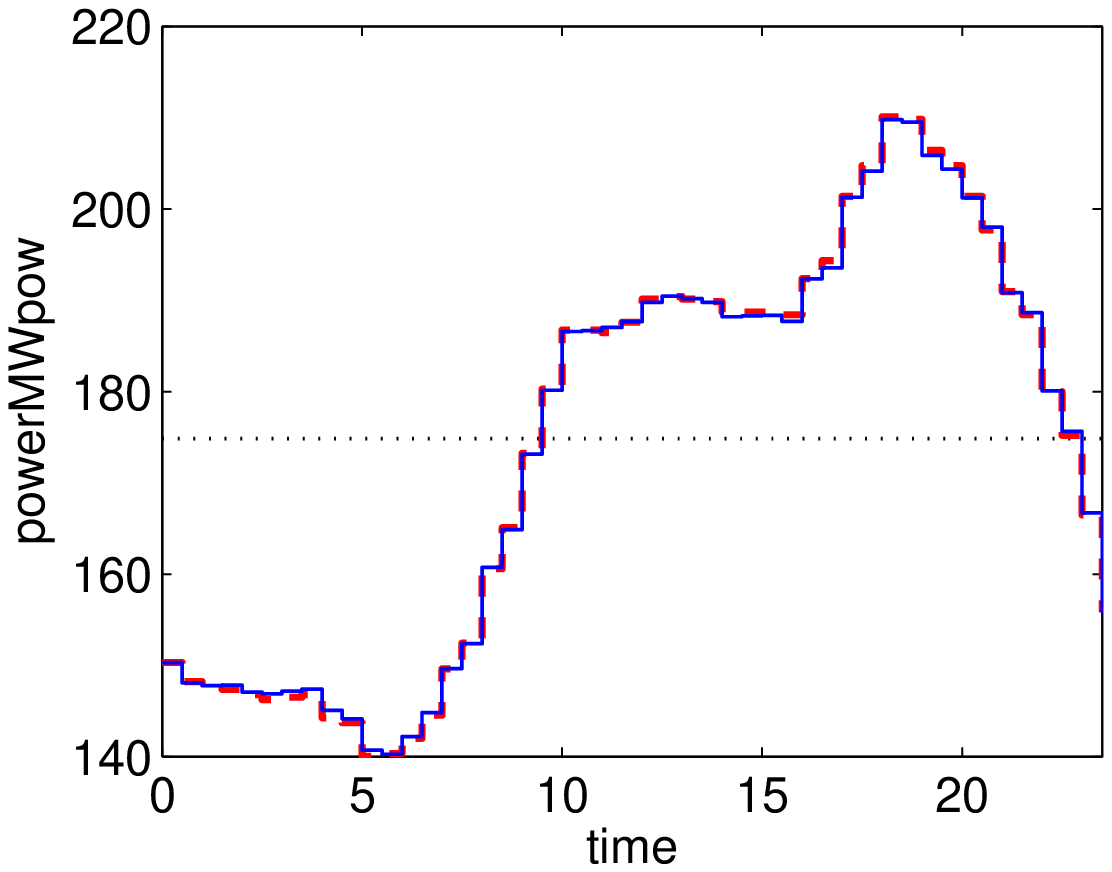}
\label{fig_tracking_loc_dyn}
}
\subfigure[DMPC and GLOBAL--REF]{
\psfrag{time}{{\footnotesize{$t$ [h]}}}
\psfrag{powerMWpow}{{\footnotesize{Power [MW]}}}
%\psfrag{aaaaaaaaaaaaaaaaaa}{{\scriptsize{Produced power}}}
%\psfrag{bbbbbbbbbbbbbbbbbb}{{\scriptsize{Reference power}}}
%\psfrag{cccccccccccccccccc}{{\scriptsize{Steady state power}}}
%\includegraphics[width=0.45\columnwidth]{power_tracking_DMPC_global_ref.eps}
\includegraphics[width=5.3cm]{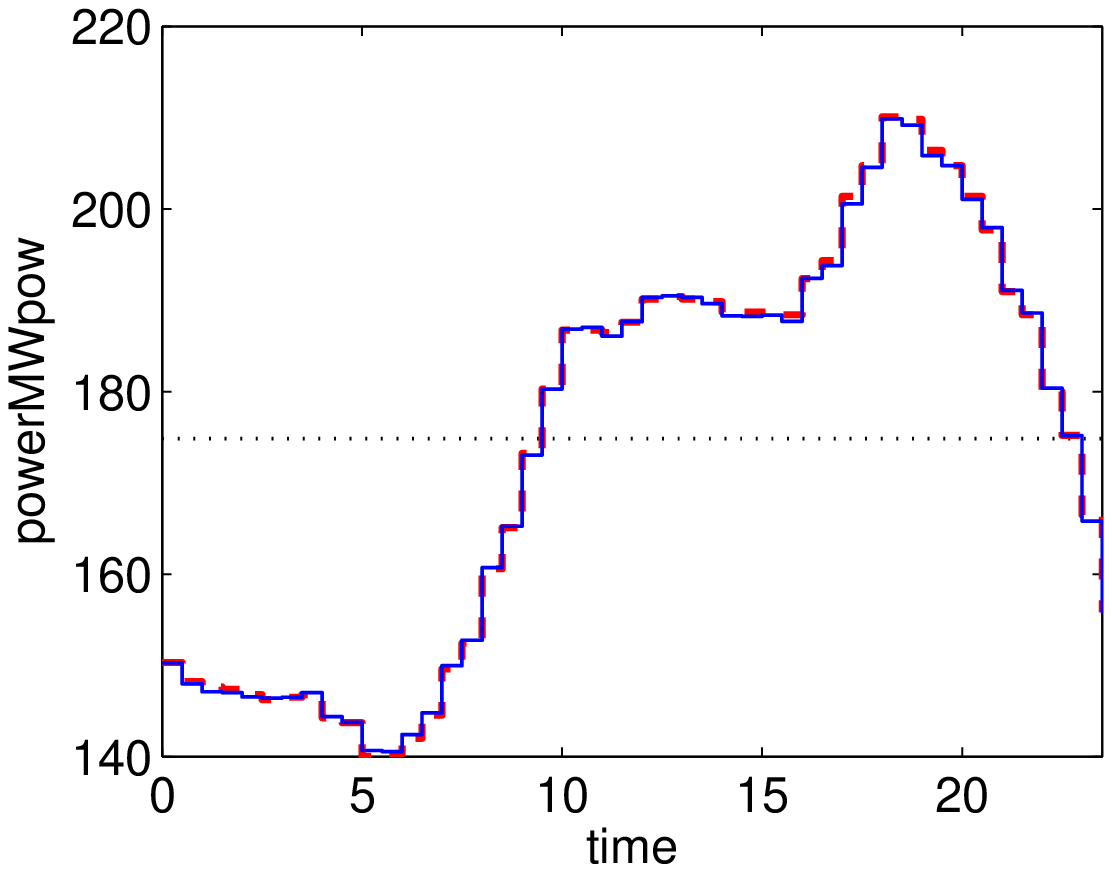}
\label{fig_tracking_global}
}
\caption{Comparison of power reference tracking performance using DMPC
  and decentralized MPC approaches. Solid lines: produced power,
  dashed lines: reference power, dotted lines: steady state power.}
\label{fig_tracking}
\end{figure}

The power reference tracking results are plotted in
Figures~\ref{fig_tracking_global}--\ref{fig_tracking_decMPC} where
the full power reference and the sum of the local power productions
are plotted. The scheme
GLOBAL--REF achieves very good tracking performance, while the scheme
LOC--REF--STAT shows a significant deterioration in tracking
performance. However, the introduction of the possibility to exchange
power references in LOC--REF--DYN between subsystems restores the very good
tracking performance while keeping the computations distributed. The
tracking performance of 
the decentralized MPC approach is very poor, due to the lack of communications.
Hence, it is recommended not to use a decentralized MPC
approach, unless communication is prohibited due to the lack of communication facilities or due to the policy of different authorities.

In \ref{app:figInputConstr} and \ref{app:figOutputConstr}
there are figures that
show the input and output evolutions and the corresponding constraints
with the scheme LOC--REF--DYN. We can observe that all constraints are
satisfied despite disturbances, model mismatch, and the use of an observer.
For the schemes GLOBAL--REF and LOC--REF--STAT all the constraints on the
inputs and outputs are also satisfied.

During the simulations, it is observed that all schemes achieve stable closed-loop behaviours, which can be explained that the HPV system is already marginally stable and does not have critical dynamics, and the prediction horizon is long enough so that the MPC controllers do not introduce instability to the closed loop. Note that neither the centralized MPC nor the distributed or decentralized MPC approaches used in this simulations employ a method that provides guaranteed stability to the closed-loop system, since this property is beyond the scope of this paper. Based on the techniques for distributing the computation and improving the efficiency of the algorithm that are proposed in this paper, one can further incorporate other MPC schemes that guarantee the closed-loop stability, which could be important for other types of applications where there are large mismatch between the nonlinear and the linearized models.

\subsection{Computational efficiency/accuracy}

In Table~\ref{tab_computation} we provide a comparison of the execution
times of the centralized MPC problems \eqref{eq:HPVopt}. We compare the distributed
Algorithm~\ref{BNB_alg} to the solver CPLEX when solving
\eqref{eq:HPVopt}, i.e., with power-division
GLOBAL--REF in Algorithm~\ref{BNB_alg}. To solve this problem using
CPLEX, an MIQP formulation is used. In every iteration of
Algorithm~\ref{BNB_alg} the relaxed problem is solved twice. We also
compare the above execution times to the case when we solve the first relaxed 
problem in Algorithm~\ref{BNB_alg}, which is a QP, using CPLEX.
At each sampling step, the same problem is solved, and
the execution time $t$ 
is measured. Although in this example the
solvers easily solve the problem within the time frame of the sampling time, we can see that
the computation time for our MATLAB-implemented
algorithm is always lower than the C-implemented CPLEX for both the
MIQP and QP cases.

\begin{table}[h!]
\centering
\caption{Comparison of computation time between
  Algorithm~\ref{BNB_alg} and CPLEX for 48 instance of the same problem}
{
\begin{tabular}{|c|c|c|c|}
\hline  & Algorithm~\ref{BNB_alg}  &  CPLEX for MIQP & CPLEX for QP\\
\hline min $t$ (s) & 0.023 &  0.087 & 0.049 \\ 
%\hline 
max $t$ (s) & 0.086 & 0.121 & 0.089 \\ 
%\hline 
average $t$ (s) & 0.054 & 0.098 & 0.063 \\
%\hline 
std dev $t$ (s) & 0.017 & 0.009 & 0.009 \\
\hline
\end{tabular}}\\
\label{tab_computation}
\end{table}
As previously discussed, Algorithm~\ref{BNB_alg} cannot
guarantee that the global optimum for \eqref{eq:HPVopt} is found. However, in the
DMPC simulations presented in this section the global optimum of
\eqref{eq:HPVopt} is found at every 
sampling step using Algorithm~\ref{BNB_alg}.

\subsection{Communication requirements}

The sizes of the optimization problems using power reference division
GLOBAL--REF, LOC--REF--DYN or LOC--REF--STAT are almost equal.
Comparing GLOBAL--REF to LOC--REF--STAT we get some additional
constraints due to the power reference division and comparing
LOC--REF--DYN to LOC--REF--STAT we get some additional decision
variables $\delta_{ij}$ to enable distributed
power reference re-assignment.

In Table~\ref{tab_communications} the number of iterations
$n_{\textrm{iter}}$ needed to
obtain the solution is presented. The average and max values of
$n_{\textrm{iter}}$ and the standard deviation are computed using 48
simulation steps, i.e., 24 hours.

\begin{table}[h!]
\centering
\caption{Number of iterations to solve the MPC optimization in one step}
{
\begin{tabular}{|c|c|c|c|c|}
\hline  & Alg.~\ref{alg:accGrad} with  &  Alg.~\ref{alg:accGrad} with  & Alg.~\ref{alg:accGrad} with  \\
  & \textbf{{GLOBAL--REF}} & \textbf{{LOC--REF--DYN}} & \textbf{{LOC--REF--STAT}} \\
%& \textbf{\scriptsize{GLOBAL--REF}} & \textbf{\scriptsize{LOC--REF--DYN}} & \textbf{\scriptsize{LOC--REF--STAT}} \\
\hline average $n_{\textrm{iter}}$  &  311.3 &  579.1 & 942.5  \\ 
%\hline 
max $n_{\textrm{iter}}$  &  498 &  1054 & 2751  \\ 
%\hline 
std dev $n_{\textrm{iter}}$  &  93.8 &  210.9 & 440.8  \\ 
\hline 
\end{tabular}}\\
\label{tab_communications}
\end{table}

We can notice that different DMPC schemes converge with different average
numbers of iterations. The reason is that for LOC--REF--STAT it is more
difficult to satisfy the different 1-norm terms with equality, i.e., to
follow the local power references. This implies that the corresponding
dual variable $\nu$ becomes large (close or equal to $\gamma$) and it
takes more iterations to achieve convergence. As a result, the scheme
LOC--REF--STAT 
with a simpler communication structure might require more
communication resources than e.g., GLOBAL--REF,
which has a more complicated communication structure but needs fewer
iterations.

In order to estimate the total time required for communications within each sampling time, we now assume the worst case happens in every iteration of Algorithm~\ref{BNB_alg}, in which Algorithm~\ref{alg:accGrad} has to be executed two times. In Algorithm~\ref{alg:accGrad}, also assume the worst case that every primal and dual variable has to be exchanged between distributed controllers, with prediction horizon $N=10$ there are $10 \times (44+65) = 1090 $ variables to be transmitted once per iteration. Let each variable be a 32-bit floating-point, then the total time it would take for transmitting exchanged variables in 1000 iterations is:
\begin{align}
 2 \times 1090 \times 32 \times 1000 = 69,760,000 (\mathrm{bits})
\end{align}
or roughly 70 Mbits. With a decent wireless network that can connect each two nodes with a rated transfer as 7 Mbps, the total time for communications is less than 10 seconds for one thousand iterations. Note that in practice, there should be more communication delays due to the initialization of transmissions. Since the communication time is considerably shorter than the sampling time of 30 minutes, the iterative methods taking about one thousand iterations sampling time can still be implemented in real time.

The scheme
LOC--REF--DYN performs very well in terms of communication, computation, as
well as performance aspects and is therefore the chosen candidate for distributed implementation for the given case study.

\section{Conclusions and future work}\label{sec_conclusions}
The proposed distributed MPC approach has been applied to the power
reference tracking problem of the HD-MPC hydro power valley benchmark.
Two distributed schemes have been compared to centralized
and decentralized MPC methods. We have provided relaxations and
approximations for the original nonlinear nonsmooth problem
formulation as well as proposed a way to follow a
centralized power reference in a distributed fashion. Furthermore, we have presented a practical
branch-and-bound algorithm that solves all optimization problems encountered in the
simulations and achieves as good performance as the centralized MPC that is known to have global optimum.
The simulation results show that the
introduced approximations and relaxations capture the behavior of
the system well and that very good control performance is achieved.
Finally, a comparison to state-of-the-art optimization software (CPLEX) shows
that the proposed algorithm has significantly better execution times in general.

As the next step before implementation in real plants,
the proposed distributed MPC approach should be tested against
different hydraulic scenarios and other HPV setups. To cope with varying water flows
entering the system, these should be estimated and compensated for.
Furthermore, a weather model could be included
that estimates the future inflows to the system.

\section{Acknowledgments}

The authors were supported by the European Union Seventh Framework STREP
project ``Hierarchical and distributed model predictive control
(HD-MPC)" with contract number INFSO-ICT-223854, the European Union
Seventh Framework Programme [FP7/2007-2013] under grant agreement no.
257462 HYCON2 Network of Excellence, the BSIK project ``Next Generation Infrastructures (NGI)", and the Swedish Research Council through the Linnaeus center LCCC.

\bibliographystyle{elsarticle-harv}
\bibliography{Doan_bib}
\newpage
\begin{appendix}

\section{}\label{app:figInputConstr}

% \begin{figure}[h!]
%  \begin{center}
%   \includegraphics[width=0.9\columnwidth]{input_constraints2.eps}
%   \caption{Input constraint satisfaction using Algorithm~\ref{BNB_alg} and
%     power division LOC--REF--DYN. Dash-dotted lines: upper bounds,
%     dashed lines: lower bounds.}
%   \label{fig_input_constraints}
%  \end{center}
% \end{figure}

\begin{figure}[h!]
 \begin{center}
\psfrag{time}{\footnotesize $t$ [h]}
\psfrag{flowflowflo0}{\footnotesize $q_{{\textrm{D6}}}$[m$^3$/s]}
\psfrag{flowflowflo1}{\footnotesize $q_{{\textrm{D5}}}$[m$^3$/s]}
\psfrag{flowflowflo2}{\footnotesize $q_{{\textrm{D4}}}$[m$^3$/s]}
\psfrag{flowflowflo3}{\footnotesize $q_{{\textrm{D3}}}$[m$^3$/s]}
\psfrag{flowflowflo4}{\footnotesize $q_{{\textrm{D2}}}$[m$^3$/s]}
\psfrag{flowflowflo5}{\footnotesize $q_{{\textrm{D1}}}$[m$^3$/s]}
\psfrag{flowflowflo6}{\footnotesize $q_{{\textrm{C2}}}$[m$^3$/s]}
\psfrag{flowflowflo7}{\footnotesize $q_{{\textrm{T2}}}$[m$^3$/s]}
\psfrag{flowflowflo8}{\footnotesize $q_{{\textrm{C1}}}$[m$^3$/s]}
\psfrag{flowflowflo9}{\footnotesize $q_{{\textrm{T1}}}$[m$^3$/s]}
  \includegraphics[width=0.9\columnwidth]{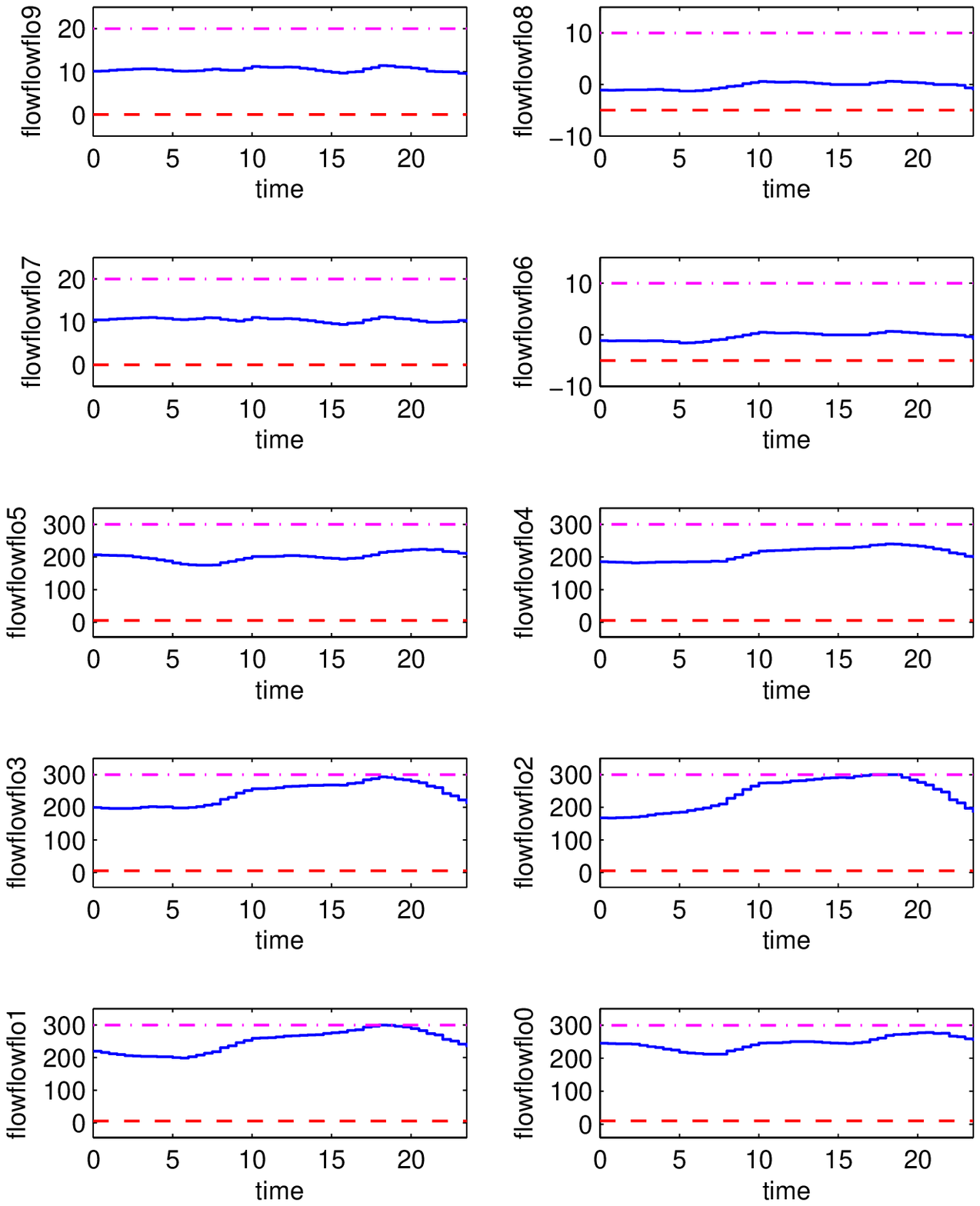}
  \caption{Input constraint satisfaction using Algorithm~\ref{BNB_alg} and
    power division LOC--REF--DYN. Dash-dotted lines: upper bounds,
    dashed lines: lower bounds.}
  \label{fig_input_constraints}
 \end{center}
\end{figure}

\newpage
\section{}\label{app:figOutputConstr}
% \begin{figure}[h!]
%  \begin{center}
%   \includegraphics[width=0.98\columnwidth]{output_constraints2.eps}
%     \caption{Output constraint satisfaction using
%       Algorithm~\ref{BNB_alg} and power division LOC--REF--DYN. Dash-dotted lines: upper bounds, dashed lines: lower bounds.}
%   \label{fig_output_constraints}
%  \end{center}
% \end{figure}

\begin{figure}[h!]
 \begin{center}
\psfrag{time}{\footnotesize $t$ [h]}
\psfrag{heighthei0}{\footnotesize $h_{{\textrm{R6}}}$ [m]}
\psfrag{heighthei1}{\footnotesize $h_{{\textrm{R5}}}$ [m]}
\psfrag{heighthei2}{\footnotesize $h_{{\textrm{R4}}}$ [m]}
\psfrag{heighthei3}{\footnotesize $h_{{\textrm{R3}}}$ [m]}
\psfrag{heighthei4}{\footnotesize $h_{{\textrm{R2}}}$ [m]}
\psfrag{heighthei5}{\footnotesize $h_{{\textrm{R1}}}$ [m]}
\psfrag{heighthei6}{\footnotesize $h_{{\textrm{L3}}}$ [m]}
\psfrag{heighthei7}{\footnotesize $h_{{\textrm{L2}}}$ [m]}
\psfrag{heighthei8}{\footnotesize $h_{{\textrm{L1}}}$ [m]}
  \includegraphics[width=0.98\columnwidth]{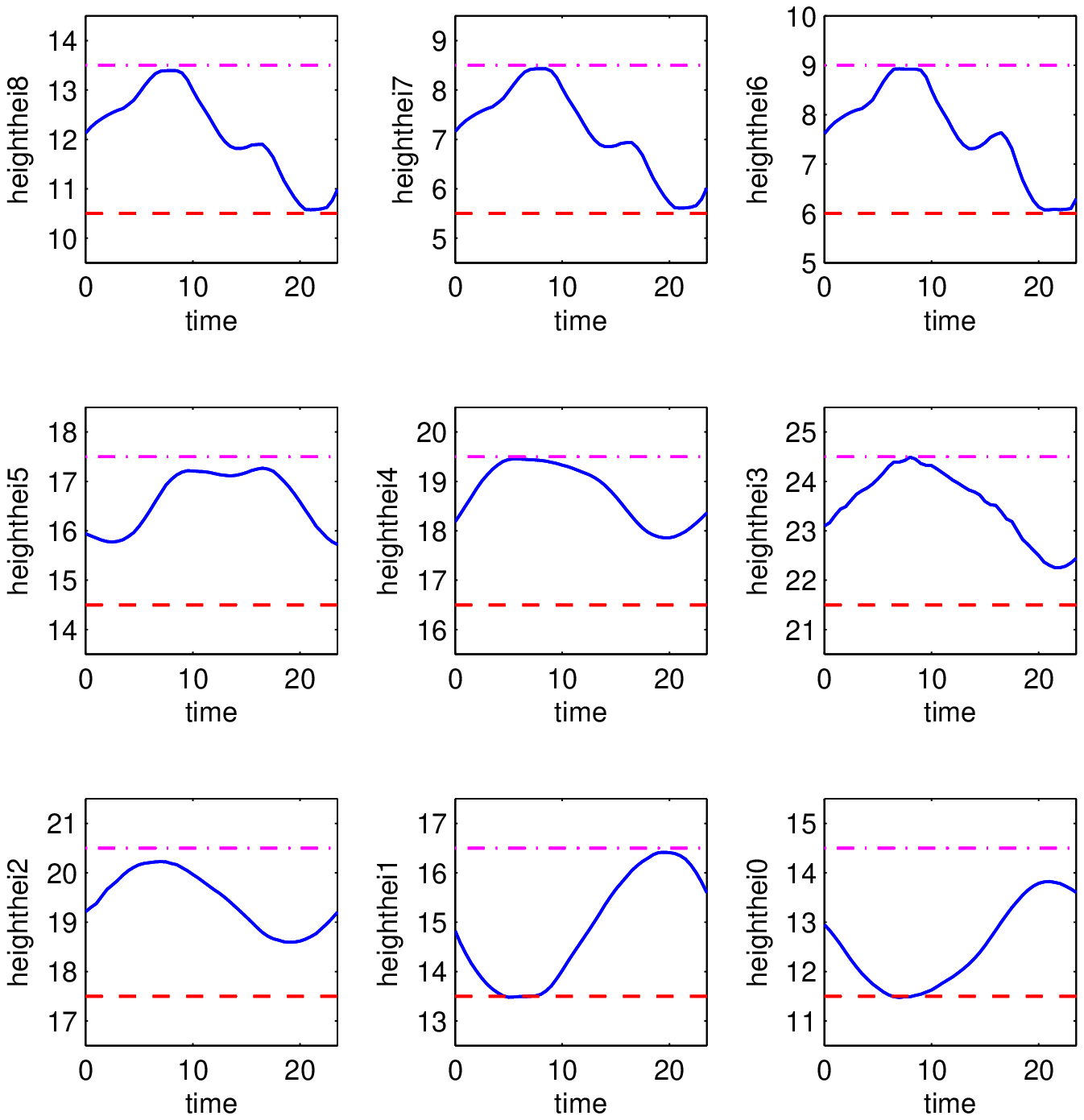}
    \caption{Output constraint satisfaction using
      Algorithm~\ref{BNB_alg} and power division LOC--REF--DYN. Dash-dotted lines: upper bounds, dashed lines: lower bounds.}
  \label{fig_output_constraints}
 \end{center}
\end{figure}

\end{appendix}

\end{document}